\documentclass[12pt]{article}
\usepackage{amsmath,amssymb,amsfonts,amsthm}
\usepackage[all]{xy}

\newcounter{item}[section]
\newcounter{kirshr}
\newcounter{kirsha}
\newcounter{kirshb}
\newenvironment{enumroman}{\setcounter{kirshr}{1}
\begin{list}{(\roman{kirshr})}{\usecounter{kirshr}} }{\end{list}}
\newenvironment{enumarab}{\setcounter{kirshb}{1}
\begin{list}{(\arabic{kirshb})}{\usecounter{kirshb}} }{\end{list}}
\newenvironment{athm}[1]{\vskip3mm\par\noindent
{\bf #1 }. \slshape }
{\upshape\par\vskip10pt minus3pt}
\newtheorem{theorem}{Theorem}[section]

\newtheorem{lemma}[theorem]{Lemma}

\newenvironment{demo}[1]{\noindent{\bf #1.}\upshape\mdseries}
{\nopagebreak{\hfill\rule{2mm}{2mm}\nopagebreak}\par\normalfont}
\theoremstyle{definition}

\newtheorem{definition}[theorem]{Definition}
\def\Q{\mathbb{Q}}
\def\C{{\mathfrak{C}}}

\def\Nr{{\mathfrak{Nr}}}

\def\Sg{{\mathfrak{Sg}}}

\def\A{{\mathfrak{A}}}
\def\B{{\mathfrak{B}}}
\def\C{{\mathfrak{C}}}

\def\CA{{\bf CA}}

\def\Dc{{\bf Dc}}

\def\RCA{{\bf RCA}}

\def\(R)RA{{\bf (R)RA}}

\def\Dc{{\bf Dc}}

\def\Q{\mathbb{Q}}
\def\Dc{{\bf Dc}}

\def\Cs{{\bf Cs}}

 \def\CA{{\sf CA}}
\def\B{{\sf B}}

\def\Nr{{\mathfrak{Nr}}}

\def\Nr{{\mathfrak{Nr}}}

\def\A{{\mathfrak{A}}}
\def\B{{\mathfrak{B}}}
\def\C{{\mathfrak{C}}}

\def\CA{{\bf CA}}

\def\RCA{{\bf RCA}}

\def\PEA{{\bf PEA}}
\def\PA{{\bf PA}}

\def\Nr{{\mathfrak{Nr}}}

\def\Sg{{\mathfrak{Sg}}}

\def\Ig{{\mathfrak{Ig}}}
\def\CA{{\bf CA}}
\def\RCA{{\bf RCA}}

\def\(R)RA{{\bf (R)RA}}

\def\Dc{{\bf Dc}}

\def\Dc{{\bf Dc}}

\def\Nr{{\mathfrak{Nr}}}

\def\Sg{{\mathfrak{Sg}}}

\def\Ig{{\mathfrak{Ig}}}
\def\CA{{\bf CA}}
\def\RCA{{\bf RCA}}

\def\(R)RA{{\bf (R)RA}}

\def\Dc{{\bf Dc}}

\def\Dc{{\bf Dc}}

\def\PA{{\bf PA}}

\def\Q{\mathbb{Q}}

\def\Q{\mathbb{Q}}
\def\Dc{{\bf Dc}}

\def\Cs{{\bf Cs}}

 \def\CA{{\sf CA}}

\def\A{{\mathfrak{A}}}
\def\B{{\mathfrak{B}}}
\def\C{{\mathfrak{C}}}

\def\Nr{{\mathfrak{Nr}}}

\def\CA{{\bf CA}}
\def\RCA{{\bf RCA}}
\def\RQEA{{\bf RQEA}}

\def\Sg{{\mathfrak Sg}}

\title{On Simple finitely generated simple polyadic finitely generated simple polyadic }
\author{Tarek Sayed Ahmed}

\begin{document}
\maketitle
 
\begin{abstract} We show that simple finitely generated polyadic equality algebras
may not be generated by a single element. We prove an analogous result for quasi-polyadic equality algebras
of infinite dimensions. In contrast, we show that any simple finitely generated 
infinite dimensional polyadic equality algeba is generated by a one element.
\footnote{ 2000 {\it Mathematics Subject Classification.} Primary 03G15.
{\it Key words}: algebraic logic, quasi-polyadic algebras, polyadic algebra} 
\end{abstract}

The two most famous algebraisations of first order logic are Tarski's cylindric algebras and Halmos' polyadic algebras.
It is commonly accepted that such algebras belong to different paradigms showing a lot of contradictory behaviour. 
In this paper, we add one such property (that mentioned in the abstract) 
to this long list establishing this dichotomy.  

Let $\alpha$ be an ordinal. Let $U$ be a set.
Let $\Gamma\subseteq \alpha,$ $E\subseteq \alpha\times \alpha$, $\tau:\alpha\to \alpha$ and $X\subseteq {}^{\alpha}U$:
$${\sf c}_{(\Gamma)}X=\{s\in {}^{\alpha}U:  \exists t\in X, t(i)=s(i) \text { for all }i\notin \Gamma\},$$
$${\sf s}_{\tau}X=\{s\in {}^{\alpha}U: s\circ \tau\in X\},$$
$${\sf d}_E=\{s\in {}^{\alpha}U: s_i=s_j \text { for } (i,j)\in E\}.$$
The operations are called {\it cylindrifications}, {\it substitutions}, and  {\it diagonal elements}, respecively.
We use simplified notation for some of these operations
${\sf c}_i={\sf c}_{(\{i\})}$ and ${\sf d}_{ij}={\sf d}_{\{(i,j)\}}$.
Let $S$ be the operation of forming subalgebras,
$P$ be that of forming products, and $H$ be the operation of forming homomorphic images.
Then 
\begin{enumroman}
\item $\RCA_{\alpha}=SP\{(\B(^{\alpha}U), {\sf c}_i, {\sf d}_{ij})_{i,j< \alpha}: U \text{ is a set }\}.$
\item $\RQEA_{\alpha}=SP\{(\B(^{\alpha}U), {\sf c}_i, {\sf d}_{ij}, {\sf s}_{[i,j]})_{i,j< \alpha}: U \text{ is a set }\}.$
\end{enumroman}

Algebras that prove handy at counterexamples are weak set algebras. These are set algebras with units of the form
$V={}^{\alpha}U^{(p)}=\{s\in {}^{\alpha}U |\{i\in \alpha: s_i\neq p_i\}|<\omega\}$ where $p\in {}^{\alpha}U$, and the operations 
are the relativization of the above defined operations on $V$.

$\Dc_{\alpha}$ is the class of dimension complemented cylindric algbras. $\A\in \Dc_{\alpha}$, if $\Delta x\neq \alpha$ for all $x\in \A$.
Here $\Delta x=\{i\in \alpha: {\sf c}_ix\neq x\}$ is called the dimension set of $x$.
It is known that, $\Dc_{\alpha}$ for infinite $\alpha$ is a nice proper generalization of locally finite algebras 
(the algebraic counterpart of first order logic, cf \cite{HMT2} section 4.3), sharing a lot of its properties \cite{s2}. In this case 
cylindric and polyadic algebras coincide. 
It is not hard to show that in this case simple finitely generated algebras are generated by a single 
element.

Andr\'eka and N\'emeti have shown that this property for cylindric algebras 
does not generalize to $\Cs_{\alpha}^{reg}$, even when $\alpha$ is finite $\geq 2$.

\section{The quasipolyadic case}

We use a fairly diect modification of their construction to prove the quasipolyadic equality case.
Let $\alpha$ be an ordinal $>0$. $\mathbb{Q}$ denotes the set of rational numbers.
Let $r\in \mathbb{Q}$, then $\bold r=\alpha\times \{r\}$. $^{\alpha}\mathbb{Q}^{(\bold r)}$ is the weak space 
$\{s\in {}^{\alpha}\mathbb{Q}: |]\{i\in \alpha: s_i\neq r\}|<\omega\}$ and 
$V={}^{\alpha}\mathbb{Q}^{(\bold 0)}$. Throughout $r\in {}^{\alpha}\mathbb{Q}$ and $t\in \mathbb{Q}$.
$$[t,r]=\{s\in V: t=\sum\{r_is_i:i<\alpha\}\},$$
$$Pol=\{[t,r]: r\in {}^{\alpha}\mathbb{Q}^{(\bold 0)}\cup {}^{\alpha}\mathbb{Q}^{(\bold 1)}\},$$
$$Po=\{[t,r]\in Pol: \{0\}\subseteq Rg(r)\}\cup \{{\sf d}_{ij}: i<j<\alpha\},$$
$$Pof=\{[t, \bold 1]: t\in \mathbb{Q}\}\text { and }\A=\Sg^{\C}Pol.$$
Her we are taking the closure under all polyadic operations. Let $Per_{\alpha}$ be the set of all transpositions $\{{[i,j]}: i,j\in \alpha\}$.
For a set $X$, let $X^S=\{{\sf s}_{\tau}x: x\in X, \tau\in Per_{\alpha}\}$.
Then clearly $X\subseteq X^S$ and $(X^S)^S=X^S$.
Then, it is absolutely straightforward to show that  $Pol^S=Pol$ and $Po^S=Po$ and $Pof^S=Pof$.
Also $a_n=[n,\bold 1]$, $1^A=[0,\bold 0]$, $0^A=[1,\bold 0]$ $Po\cup Pof\subseteq Pol\subseteq \wp(V)=C$,
$\{a_n: n\in \omega\}\subseteq Pof $ and $\{0,1\}\subseteq Pol \sim (Po\cup Pof)$.
Let $X\subseteq C$ be arbitary, then
$$X^*=\{\prod \{p_i:i<n\}: n\in \omega \text { and }(\forall i <n)(p_i\in X\text { or } -p_i\in X)\},$$
$$X^{**}=\{\sum\{\pi_i:i<n\}: n\in \omega \text { and } (\forall i<n)\pi_i\in X^*\},$$
$$G(X)=(X^S\cup Po)^{**}.$$
Note the definition of $G(X)$ is different than that in \cite{AN}, it is in fact bigger.
In what follows, we show that the proof in \cite{AN} survives introducing 
substitutions coresponding to transpositions. 

\begin{lemma} $(\forall Y\subseteq G(Pol))(G(Y)\subseteq \C$
\end{lemma} 
\begin{demo}{Proof} \cite{AN} Lemma 1.
The two technical lemmas 1.1 in  and 1.2 in \cite{AN} are the same. 
$G(Y)$ is closed all the cylindric operations is done exactly like the proof for the $\CA$ case.
We need to check that $G(Y)$ is closed under substitutions.
Let $p\in G(Y)$. Then $p\in (Y^S\cup Po)^{**}$. Since $Y^S$ and 
$Po$ are closed under substitutions, and substitutions are Boolean endomorphisms, then we are done. 
\end{demo}
\begin{lemma} $\A$ is simple, and every subalgebra of $\A$ is simple
\end{lemma}
\begin{demo}{Proof} Like the proof of claim 1.1 p. 868 in \cite{AN}. 
\end{demo}

\begin{lemma} 
\begin{enumroman}
\item $nat:Pof^{**}\to G(Pof)/Poz$ is an isomorphism of Boolean algebras, that is $Poz\cap Pof^{**}=\{0\}$
\item $Pof^{**}$ is an atomic Boolean algebra, and the set of its atoms is $Pof$.
\item Let $Y\subseteq Pof$. Then $G(Y)\cap Pof=Y^S$ and $G(Y)\cap Pof^{**}=(Y^S)^{**}.$
\end{enumroman}
\end{lemma}
\begin{demo}{Proof} We need to check only the last statement, which is slightly different from the corresponding one in \cite{AN}, 
namely that formulated in lemma 2 (iii) p. 868 of \cite{AN}, for we are allowing substitutions. 
Let $Y\subseteq Pol$. First we show that $G(Y)\cap Pof^{**}\subseteq (Y^S)^{**}.$
Let $d\in G(Y)\cap Pof^{**}$. Then $Y^S\cup P$ generates $d$ in $G(Pof)$. Consider the factor algebra $G(Pof)/Poz$.
Then $Y^S\cup P$ generates $d$ in $G(Pof)$. Then $Y/Poz$ generates $d/Poz$ in $G(Pof)/Poz$. 
But $nat:Pof^{**}\to G(Pof)/Poz$ is an isomorphism. Since $Y^S\cup \{d\}\subseteq Pof^{**}$, this
implies that $Y^S$ generates $d$ in $Pof^{**}$, i.e $d\in Y^S$. 
We have proved that $G(Y)\cap Pof^{**}=Y^S$. Let $y\in G(Y)\cap Pof$. Then $y\in (Y^S)^{**}$. But $Pof$ is the set of atoms of
$Pof^{**}$. Since $|Pof|\geq \omega$, 
no set of generators generates new atoms in $Pof^{**}$. Hence $y\in  Y^S$ We have seen that $G(Y)\cap Pof=Y^S$.
\end{demo}

\begin{definition} Let $Neg=G(0)\sim Poz$.
\end{definition}
\begin{athm}{Claim 4} 
\begin{enumroman}
\item $Poz\subseteq G(0)$
\item $Neg$ is an ultrafilter of the Boolean algebra $(G(0), \cap, -)$. $Neg$ is the filter generated by $\{-p: p\in Po\}$.
\item $(\forall g\in Pof^{**}\sim \{0\})(\forall \sigma\in G(0)[g\cdot \sigma\in Poz\Leftrightarrow \sigma \in Poz].$
\item $(\forall Y\subseteq Pof)
(\forall g\in G(Y))(\exists \rho\in G(0)(\exists n\in \omega)(\exists v\in {}^nNeg)(\exists y\in {}^nY)\rho+\sum_{i\in n}v_iy_i\in \{g,-g\}.$
\end{enumroman}
\end{athm}
\begin{demo}{Proof} \cite{AN} Lemma 3. We need to check the last item only. Let $Y\subseteq Pof$ and $g\in G(Y)=(Y\cup Po)^{**}$. 
Then there is a finite $W\subseteq Y^S$ such that $g\in G(W)$ and $(\forall W_0\subseteq W)(g\notin G(W_0)$. Let $n=|W|$ and
$\{y_i: i\in n\}$. Then every $y_i$ is either an element of $Y$ or a substitution corresponding to a transposition 
of an element in $Y$ and they are all distinct. 
Since $y_i\neq y_j\implies y_i\cap y_j=0$, because $Y^S\subseteq Pof$, and the latter is an antichain. Since $G(W)=(W^S\cup Po)^{**}$, we have that
$$g=\sum_{i\in n}y_i\cdot \sigma_i+\sigma_n\cdot \prod_{i\in n}-y_i$$
for some $\{\sigma_i: i\in n\}\subseteq G(0)=Po^{**}$.
Then $$-g=\sum_{i\in n} y_i\cdot -\sigma_i + -\sigma_n\cdot \prod\{-y_i: i\in n\}.$$
Assume $\sigma_n\in Poz$. Then $\sigma_n.\prod\{-y_i:i\in n\}\in Poz$ and
$$(\{\sigma_i:i\in n\}\cap Poz\neq 0\implies \{y_i\cdot \sigma_i: i\in n\}\cap Poz\neq 0)$$
by $W^S\subseteq Pof$. Then $(\forall W_0\subset W)(g\notin G(W_0)$ and $Poz\subseteq G(0)$ imply $\{\sigma_i:i\in n\}\subseteq Neg$
so we get the desired form. If $\sigma_n\notin Poz$, then $-\sigma_n\in Poz$ and we can work analogously with $-g$.
\end{demo}
Let $$\A_n=\Sg^{\C}\{a_i:i< 2^n\}.$$
That is $\A_n$ is the quasi-polyadic equality subalgebra of $\C$, genertaed by the $a_i$'s $i<2^n$.
Recall that $\A_n$ is simple. Then we have our first main Theorem:
\begin{athm}{Theorem 2} $\A_n$ cannot be generated by $n$ elements, but can be 
generated by $n+1$ elements
\end{athm}

\begin{demo}{Proof} \cite{AN} Claim 1.2 p. 871.
Fix $n$, and let $Y=\{a_i: i< {}^n2\}$. First we show that
$$(\forall X\subseteq A_n)[|X|\leq n\implies \Sg^{\C}X\neq A_n]$$
Suppose to the contrary  that $\Sg^{\C}X=\A_n$ for some $X=\{g_i:i<n\}$. Then $A_n\subseteq G(Y)$ and $A_n\subseteq G(X)$.
Let $\{g_i:i\in  n\}=X$ and let $j<n.$
Then we amy suppose that
$$(\exists I_j\subseteq 2^n)(\exists \rho_j\in G(0))(\exists v^j:I_j\to Neg)g_j=\rho_j+\sum\{v_i^ja_i: i\in I_j\}.$$
Let these $I_j$ $\rho_j$ and $v^j$ be fixed for every $j\in n$. By $Y\subseteq G(X)\subseteq G(\{a_i: i\in \bigcup\{I_j:j\ in n\}$ we have that 
$\bigcup\{I_j: j\in n\}=2^n$.
Then 
$$\exists k, h\in 2^n(\forall j\in n)[[k\in I_j\text { iff }h\in I_j\}$$
by the same reasoning in \cite{AN} p. 872.
Let those $k$ and $h$ be fixed. By $a_k\in G(X)$ we have that $a_k=\sum B'$ for some finite $B'\subseteq (X^S\cup Po)^*$.
Let $B=[B'\cap (X\cup Po)^{**}]$. Then, like the proof in \cite{AN} p. 872, it can be shown that $B\subseteq G(0)$. 
Thus $B'\subseteq G(0)$ since the latter is closed under substitutions corresponding to
transpositions.
Thus $a_k\in G(0)$ which is a contradiction. The rest of the proof is essentially the same as that in \cite{AN} p.872.
\end{demo}

Finally, we should mention that the above proof survives if we replace $\Q$ by any field of characteristic $0$.
However the proof does not work when $\Q$ is replaced by a finite field. 

Now we turn to the polyadic paradigm. $\PA_{\alpha}$ denotes the class of polyadic algebras of dimension $\alpha$, while $\PEA_{\alpha}$
denotes the class of polyadic equality algebras of dimension $\alpha$. From now on $\alpha$ will be only infinite.
We follow \cite{HMT2} for the abstract axiomatixation of such algebras.

\begin{athm}{Definition 4} Let $\A\in \PA_{\alpha}$. Let $J\subseteq \alpha$, an element $a\in A$ is {\it independent of $J$} if ${\sf c}_{(J)}a=a$.
$J$ supports $a$ if $a$ is independent of $\alpha\sim J$. We write $\Delta a$ for the least $J$ that supports $a$; 
$\Delta a$ is called the dimension of $a$.
\end{athm}

\begin{athm}{Definition 5} Let $J\subseteq \beta$ and 
$\A=\langle A,+,\cdot ,-, 0, 1,{\sf c}_{(\Gamma)}, {\sf s}_{\tau}\rangle_{\Gamma\subseteq \beta ,\tau\in {}^{\beta}\beta}$
be a $\PA_{\beta}$.
Let $Nr_J\B=\{a\in A: {\sf c}_{(\beta\sim J)}a=a\}$. Then
$$\Nr_J\B=\langle Nr_{J}\B, +, \cdot, -, {\sf c}_{(\Gamma)}, {\sf s}'_{\tau}\rangle_{\Gamma\subseteq J, \tau\in {}^{\alpha}\alpha}$$
where ${\sf s}'_{\tau}={\sf s}_{\bar{\tau}}.$
The structure $\Nr_J\B$ is an algebra, called the {\it $J$ compression} of $\B$.
When $J=\alpha$, $\alpha$ an ordinal, then $\Nr_{\alpha}\B\in \PA_{\alpha}$ and is called the {\it neat $\alpha$ reduct} of $\B$ 
and its elements are called 
$\alpha$-dimensional.
\end{athm}

\begin{athm}{Lemma}
\begin{enumarab}
\item Let $\A\in \PEA_{\alpha}$. Then for all $\beta>\alpha$ there exists $\B\in \PEA_{\beta}$ such that $\A\subseteq \Nr_{\alpha}\B$.
Furthermore, every element of $\B$ is of the form ${\sf s}_{\sigma}^{\B}a$ where $a\in A$ , $\sigma\in {}^{\beta}\beta$ and $\sigma\upharpoonright \alpha$ 
is one to one.
\item If $A$ is a generates $\B$ and $\A$ is simple, then so is $\B$
\end{enumarab} 
\end{athm} 

\begin{demo}{Proof} For a subset $X$ of an algebra $\C$, we write $\Ig^{\A}X$, for the ideal generated by $X$.
We show that if $I$ is an ideal in $\B$, then the ideal generated by $I\cap A$ in $\B$ coincides with $I$. We have 
$\A=\Nr_{\alpha}\B$. Only one incusion is non-trivial. Let $x\in I$.
The ${\sf c}_{(\Delta x\sim \alpha)}^{\B}x\in A$, hence in $I\cap A$. 
Since $x\leq {\sf c}_{(\Delta x\sim \alpha)}^{\B}x$, we have $x\in \Ig^{\B}(I\cap A)$. 
So it follows that if $I$ is a ideal in $\B$, and $\A$ is simple, so that $\A\cap I=\{0\},$ then
$I=\{0\}.$
\end{demo}

\begin{athm}{Theorem} Every simple $\PEA_{\alpha}$ generated by finitely many elements is generated by a single element
\end{athm} 
\begin{demo}{Proof} Let $\A\in \PEA_{\alpha}$ be simple, and assume that it is  generated by two elemens $x$ and $y$.
Let $\beta$ be an ordinal such that $\beta$ contains at least two distinct elements $k,l$ not $\in \alpha$, and 
let $\B$ be a dilation of $\A$ so that that $\A=\Nr_{\alpha}\B$, and $A$ generates $\B$, so that $\B$ is, 
in fact,  generated by $x$ and $y$. 
Let $k,l\in \beta\sim \alpha$. Let $b= x\cdot {\sf d}_{kl}+y\cdot {\sf d}_{kl}.$ Then we claim that
$\B=\Sg^{\B}{b}$. First not that $b\cdot {\sf d}_{kl}=x\cdot {\sf d}_{kl}$ and $z-\cdot {\sf d}_{kl}=
y\cdot -{\sf d}_{kl}.$
Then $${\sf c}_k(z\cdot {\sf d}_{kl})={\sf c}_k(x\cdot {\sf d}_{kl})=x\cdot {\sf c}_k{\sf d}_{kl}=x,$$
and
$${\sf c}_k(z\cdot -{\sf d}_{kl})={\sf c}_k(y\cdot -{\sf d}_{kl})= y\cdot {\sf c}_k-{\sf d}_{kl}=y\cdot {\sf c}_0\cdot -{\sf d}_{01}=y.$$
Hence $\B$ is generated by $b$. Since $\B$ is a dilation of $\A$, then there exists $\sigma\in {}^{\beta}\beta$, $\sigma$ is 
one to one on $\alpha$ and  $a\in A$ such that 
${\sf s}_{\sigma}^{\B}a=b$. We claim that $a$ generates $\A$. 
Let $\Omega=\beta\sim \alpha$. It suffices to show that if $x\in \B$, ${\sf c}_{(\Omega)}^{\B}x=x,$ then $x\in \Sg^{\A}\{a\}$. 

We first treat several cases separately. Call a term $\tau$ positive if it is 
built up of succesive applications of substitutions and cylindrifications  
and has no occurences of complementation.

We start studying some positive terms applied $b$, and show  that when the output is in $\A$, then 
it is actually in $\Sg^{\A}\{a\}.$

We proceed by induction on the number of unary operations. We start by $n=1$. 

(1) Let $x={\sf s}_{\sigma'}^{\B}b$, where $\Delta x\subseteq \alpha$. Then $x= {\sf s}_{\sigma'}^{\B}{\sf s}_{\sigma}^{\B}a$. 
Let $\mu=\sigma\circ \sigma'$. Choose $\tau\in {}^{\beta}\beta$, such that $\tau\circ \mu(\alpha)\subseteq \alpha$, 
and $\tau\upharpoonright \alpha\subseteq Id.$
Then since $\Delta x\subseteq \alpha$,
we get
$$x={\sf s}_{\mu}^{\B}a={\sf s}_{\tau\circ \mu}^{\B}a={\sf s}_{(\tau\circ \mu)\upharpoonright \alpha)}^{\B}a=
{\sf s}_{\tau\circ \mu}^{\A}a\in \Sg^{\A}\{a\}.$$ 

(2) Assume that $x={\sf c}_{(\Delta)}b$ where $\Delta\subseteq \beta$. Then $x={\sf c}_{(\Delta)}^{\B}{\sf s}_{\sigma}^{\B}a$. 

Choose $\mu$ and $\tau$ permutations of $\beta$ such that 
$(\mu\circ \tau)| \alpha\subseteq Id$, $\mu\circ \tau\circ \sigma(\alpha)\subseteq \alpha$ and 
$\mu\circ \tau(\Delta)\subseteq \alpha$.
Then 
\begin{equation*}
\begin{split}
{\sf c}_{(\Delta)}^{\B}{\sf s}_{\sigma}^{\A}a
&={\sf s}_{Id}{\sf c}_{(\Delta)}^{\B}{\sf s}_{\sigma}^{\B}a\\
&={\sf s}_{\mu\circ \tau)\upharpoonright \alpha}^{\A} {\sf c}_{(\Delta)}^{\B}{\sf s}_{\sigma}^{\B}a\\
&={\sf s}_{\mu}^{\B}{\sf s}_{\tau}{\sf c}_{(\Delta)}^{\B}{\sf s}_{\sigma}^{\A}a\\
&={\sf s}_{\mu}{\sf c}_{\Gamma}{\sf s}_{\tau}{\sf s}_{\sigma}a\\ 
&={\sf c}_{\Delta'}{\sf s}_{\mu}{\sf s}_{\tau}{\sf s}_{\sigma}a\\
&={\sf c}_{(\Delta')}{\sf s}_{\mu\circ \tau\circ \sigma}a\\
&={\sf c}_{(\Delta')}{\sf s}_{\mu\circ \tau\circ \sigma}a\\
\end{split}
\end{equation*}

Now we have shown that if we have a term $\tau$ of the form $f\circ f_2\circ \ldots f_n$ 
where the $f_i$'s are either  cylindrifications or substitutionsand $\tau(b)\in \A$, then $\tau(b)\in \Sg^{\A}(a)$.

Now we consider the Boolean join $+$ and complementation.
We start by complementation applied to unary terms. If $\tau$ is a unary term consisting of only substitutions, then $-\tau(b)=\tau(-b)$ 
since the substitutions are  
Boolean endomorphisms, and we are done. Now let ${\sf c}_{(\Delta)}^{\partial}$ 
be the dual of ${\sf c}_{(\Delta)}$ defined by $-{\sf c}_{\Delta}-$ 
If complementation is applied to a unary term containg cylindrifications, 
then by noting that 
$-{\sf c}_{(\Gamma)}a={\sf c}_{(\Gamma)}^{\partial}-x$, and the polyadic axioms P11, and P12 
involving the interaction of substitutions and 
cylindrifications,  hold for the duals of cylindrifications, we are done in this case, too.

Now we consider the Boolean join. Let $x\in \A$ and assume that $x=x_1+x_2$. 
Let $\Gamma=\beta\sim \alpha$. Then ${\sf c}_{(\Gamma)}x={\sf c}_{(\Gamma)}x_1+{\sf c}_{(\Gamma)}x_2$, 
hence ${\sf c}_{(\Gamma}x_i=x_i$,
so that $x_i\in A$ for $i=1,2.$ 
Now $x_i\in \B$ so we can assume that there exists unary terms (possibly) with negations, such that 
$\tau_1$ and $\tau_2$ such that $x_i=\tau_i(b).$
But then $x_i= \tau_i'(a)$ and we are done.
\end{demo}

\end{document}